\documentclass[11pt,draft,a4paper]{article}
\usepackage{amssymb}
\usepackage{amsmath}
\newtheorem{defi}{Definition}[section] 
\newtheorem{pro}[defi]{Proposition} 
\newtheorem{theo}[defi]{Theorem}
\newtheorem{coro}[defi]{Corollary}
\newtheorem{rela}[defi]{Relation}
\newtheorem{lem}[defi]{Lemma}
 
\newcommand{\mrm}{\mathrm} 
\newcommand{\ff}{\mathcal{F}} 
\newcommand{\mca}{\mathcal} 
\newcommand{\sui}[1]{\left(#1 _n\right)_{n\in \mathbb N}} 
\newcommand{\suit}[1]{\left(#1\right)_{n\in \mathbb N}} 
\newcommand{\matg}[4]{\begin{pmatrix}#1_n&#2_n\\#3_n&#4_n\end{pmatrix}}

\newcommand{\nN}{n\in \mathbb N}

\newcommand{\istb}[1]{AS(#1_n^{-1})}

\newcommand{\ninf}{\underset{n\rightarrow\infty}{\longrightarrow}} 
\newcommand{\st}[1]{{AS}^{\perp}(#1)} 
\newcommand{\ist}[1]{{AS}^{\perp}(#1^{-1})}

\newcommand{\diff}{diffeo\-morphism}
\newcommand{\Dif}{Diff$(M)$}

\newcommand{\e}{\mbox{\rm\Large e}}
\newcommand{\hop}{\vskip .3cm\noindent}
\title{Dynamical properties of the space of Lorentzian metrics.}
\author{Pierre Mounoud}
\date{}
\begin{document}
\maketitle
{\it\`A la m\'emoire de mon p\`ere.}
\begin{abstract}
We study the mechanisms of the non properness of the action of 
the group of diffeomorphisms on the space of Lorentzian metrics of a compact manifold.
In particular, we prove that nonproperness entails the presence of 
lightlike geodesic foliations of codimension $1$.
On the $2$-torus, we prove that a metric with constant curvature along one of its lightlike foliation is actually flat. This allows us to show that
the restriction of the action to the set of non-flat metrics is proper and that on the set of flat metrics of volume $1$ the action is ergodic.
Finally, we show that, contrarily to the Riemannian case, the space of metrics without isometries is not always open. 
\end{abstract}
\hrulefill
\vskip 0.1cm\noindent
\hrulefill
\section*{Introduction}
We study  the action of the group of diffeomorphisms of a compact
manifold $M$, denoted by \Dif, on the space of Lorentzian metrics of
$M$, denoted by $\mca L(M)$.  An important motivation is to understand the space of Lorentzian structures, ie the quotient space. In the Riemannian case this action is well understood, thanks to the work of 
D.G. Ebin \cite {ebin}. He proved  in particular that it is proper.

Even though properness  has no reason to hold in the pseudo-Riemannian case, the situation is rather intricate. On one hand we can single out manifolds for which this action is indeed proper. On the other hand non properness turns out to be a very special property. This is related to  the works of R. Zimmer \cite{zim}, 
G. D'Ambra and M. Gromov \cite{DGr}, S. Adams and G. Stuck \cite{AS}, or A. Zeghib \cite{zeho} and \cite{zeid}, about non compact isometry groups of closed Lorentzian manifolds.

Our first result is a generalization of the
main theorem of the article \cite{zeghib} by  A. Zeghib which 
deals with non-equicontinuous sequences of isometries of a compact Lorentzian manifold. This result enables us to know how the properness fails.
As we do not want to be too technical in this introduction, we formulate a weak version of theorem \ref{zorro} as follows.
\vskip 0.3cm\noindent
{\bf Theorem} {\sl Let $K$ and $L$ be two compact subsets of $\mca L(M)$ (for the $C^2$-topology) such that
 $\{\varphi\in \mathrm{Diff}(M) | \varphi^*K\cap L\neq \emptyset\}$ is not compact. Then both $K$ and $L$ contain metrics which
possess a Lipschitz, codimension $1$, totally geodesic and lightlike foliation.}  
\vskip .3cm
Recall that a submanifold of $(M,h)$ is lightlike if it is everywhere tangent to the light cone, ie if the restriction of the metric is everywhere degenerate.
 As in  Zeghib's theorem, the foliation has a dynamical meaning : its tangent is the set of 
 ``approximately stable vectors'' relative to a given sequence of diffeomorphisms. When $K$ and $L$ are reduced to a point, this theorem describes metrics with non compact isometry group, it is Zeghib's theorem.

 The third section of this article is devoted to compact Lorentzian surfaces, ie to the torus and the Klein bottle (since the Euler characteristic must vanish). In this case we are able to be more precise in the description of the non-properness. It comes from the fact that the metrics described by the previous theorem have constant curvature along one of their lightlike  foliations. 
Studying those metrics we obtain  the following global geometric property of Lorentzian tori.
\vskip 0.3cm\noindent
 {\bf Theorem \ref{bonne}} {\sl
A Lorentzian metric on the torus $\mathbf T^2$ whose curvature is
constant along one of its lightlike foliations is flat.} 
\vskip .3cm
Hence, on the torus the non-properness is localized on flat metrics. 
We still have to understand the restriction of the action to the set of flat metrics, that we will denote by $\mca F(\mathbf T^2)$. The volume functional is clearly invariant under the action of Diff$(\mathbf T^2)$.
 It is therefore natural to restrict the action to the set of flat metrics of volume $1$, that we denote by  $\mca F_1(\mathbf T^2)$. 
Next, in order to obtain a finite dimensional problem, we take the quotient of this set by the proper action of Diff$_0(\mathbf T^2)$, the group of diffeomorphisms isotopic to the identity. 
The quotient space is diffeomorphic to the set of quadratic forms of signature $(1,1)$ and of determinant $1$, ie to $SL(2,\mathbb R)/SO(1,1)$. Furthermore, it is well known that Diff$^+(\mathbf T^2)/$Diff$_0(\mathbf T^2)$ is isomorphic to $SL(2,\mathbb Z)$.
Therefore the problem is now to understand the action of $SL(2,\mathbb Z)$ on $SL(2,\mathbb R)/SO(1,1)$.
This action being
ergodic, we can state the following (cf. corollary \ref{plat} and theorem \ref{ergod}).
\vskip .3cm
\noindent
{\bf Theorem}\indent
{\sl {\it 1.} The action of  {\rm Diff}$(\mathbf T^2)$ on $\mca L(\mathbf T^2)\setminus\mca F(\mathbf T^2)$
is $C^2$-proper.}
\vskip .2cm
{\it 2.}
{\sl The action of {\rm Diff}$(\mathbf T^2)$ on  $\mca F_1(\mathbf T^2)$ is ergodic.}
\hop
We use the word ergodic  because, even if no measure is involved in this infinite dimensional context, it is nevertheless a good description of the situation.
\vskip .2cm
The last part is devoted to the set $\mca L_T(M)$ of metrics
with trivial isometry group. This  set is dense in $\mca L(M)$ and we want to know 
whether it is an  open set. 
Indeed, this is the case when Diff$(M)$ acts properly.
 Typical examples are the Klein bottle (cf. corollary \ref{bouteille}) or a $3$-manifold not covered by $\mathbb R^3$ (cf. corollary \ref{coro}). 
Although Diff$(\mathbf T^2)$ does not act properly, $\mca L_T(\mathbf T^2)$ is also open in $\mca L(\mathbf T^2)$.
However, in dimension greater than two, $\mca L_T(M)$ may not be open in $\mca L(M)$. A typical example is given by compact quotients of $PSL(2,\mathbb R)$.
\vskip 0.5cm
This article is a part of my doctoral dissertation (\cite{these}) written under the
supervision of 
Jacques Lafontaine. I wish to thank him for his attention throughout this work. 
I also wish to thank  Abdelghani Zeghib, for the discussions we had and his 
suggestions, and Yves Carri\`ere who kindly fixed the proof of lemma 2.2.
I thank also the referees for their constructive remarks.
\section{The Riemannian situation.}
For any (pseudo) metric $g$, the isometry group Is$(M,g)$ is a Lie group (Myers-Steenrood). It is compact when $M$ is compact and $g$ Riemannian: indeed the set of isometries is equicontinuous.
With the same kind of arguments, it can be proved that the natural action of Diff$(M)$ on the space of Riemannian metrics is proper. In particular, the orbits are closed.

A deeper result of D. Ebin is the slice theorem. It roughly says that any orbit of this action has a Diff$(M)$-invariant tubular neighborhood. An important consequence is the (lower) semi-continuity of the cardinal and the dimension of Is$(M,g)$ with respect to $g$ (cf. theorem 8.1 of \cite{ebin}). In particular the set of  metrics without isometries is open.

We shall study all these questions: properness of the action of Diff$(M)$, behavior of the orbits, properties of the set of isometry-free metrics, in the Lorentzian case, and obtain  properties very different to those in the Riemannian case.
\section{Approximately stable foliations.} 
The action of \Dif\ on $\mca L(M)$ is proper (for a given topology) if  for any compact subsets $K$ and  $L$ of $\mca L(M)$ the set $\{\varphi\in \mathrm{Diff}(M) | \varphi^*K\cap L\neq\emptyset\}$ is compact. In fact, as it is  more convenient  to work with sequences rather than subsets, we will characterize the \emph{non}-properness, for the $C^k$-topology, by the existence of two $C^k$-convergent sequences of metrics,  $\sui{h}$ and  $\sui{g}$, and
of a non-equicontinuous (ie not lying in a compact subset) sequence of \diff s $\sui{\varphi}$ such that $\varphi_n^*h_n=g_n$. We will call such a triple of sequences a  \emph{$C^k$
approximately isometric system}. We can view the approximately isometric systems
as generalizations of the non-equicontinuous sequences of isometries. To understand the non-properness of the action we study the dynamic  of those systems. Let us start with the definition (see \cite{zeghib}) of an 
approximately stable vector (relative to a sequence of diffeomorphisms), which will be the key notion of 
this section.
\begin{defi}\label{def}
Let $\sui{\varphi}$ be a sequence of \diff s.
\begin{enumerate} 
\item[$-$]
A vector  $v\in T_xM$ is said to be approximately stable if there exists a
convergent sequence  
 $\sui{v}$ with limit $v$, such that the sequence
$\suit{D\varphi_n(v_n)}$ is bounded. We will denote by
${AS}(\varphi_n)$  the set  of approximately stable vectors of
$\sui{\varphi}$ and by ${AS}(x,(\varphi_n))$ its  intersection with  $T_xM$.\\ 
We say that $v$ is strongly approximately stable  if in addition
$D\varphi_n(v_n)\rightarrow 0$.
\item [$-$]
A vector  $v\in T_xM$ is called punctually approximately stable if
in addition the sequence  
 $\sui{v}$ stays in the same tangent vector space  $T_xM$. We will denote by
${PAS}(\varphi_n)$  the set  of punctually approximately stable vectors of
$\sui{\varphi}$ and by ${PAS}(x,(\varphi_n))$ its  intersection with  $T_xM$.
\end{enumerate} 
\end{defi}
The goal of this section is to prove the following theorem. It is a
generalization of the main theorem of the article \cite {zeghib} of
A. Zeghib, which describes the dynamical
properties of the non-equicontinuous sequences of
Lorentzian isometries. 
If we compare the statements we see that the dynamical properties are mainly the same.
 It 
gives us a  ``localization'' of the non-properness.
\begin{theo}\label{zorro} 
Let $(h_n,g_n,\varphi_n)$ be a  $C^k$
approximately isometric system (with $k\geq 2$) on a compact manifold $M$ and let
$g_\infty$ be the limit of  $\sui{g}$.
Then  there exists a subsequence  $\sui{\phi}$ of $\sui{\varphi}$ such that the set
of approximately stable vectors ${AS}(\phi_n)$ is the tangent bundle
of a codimension 1 Lipschitz foliation, called the approximative stable
foliation of $\sui{\phi}$. The leaves of the foliation are geodesic and
lightlike  relatively to $g_\infty$.\\\indent
Replacing $\sui{\phi}$ by a subsequence if necessary, we can assume that the same is true for $\sui{\phi^{-1}}$ 
(the foliation ${AS}(\phi_n^{-1})$ will be  geodesic and lightlike relatively to $h_\infty$). In this case, 
if $v \in TM \setminus
 {AS}(\phi_n)$, then $D\phi_n v$ tends to $\infty$, and converges projectively to 
${AS}^\perp(\phi_n^{-1})$. 
The convergence 
is uniform on compact subsets of $TM \setminus {AS}(\phi_n)$. 
\\\indent
Moreover, if $\sui{\sigma}$ and  $\sui{\sigma'}$ are convergent sequences of
functions on $M$ whose limits are respectively  $\sigma$ and $\sigma'$
 and such that $\sigma'_n=\sigma_n\circ \varphi_n$, 
then $\sigma'$ is constant along $\st{\phi_n}$ (and $\sigma$ along $AS^\perp (\phi_n^{-1})$).
\end{theo} 
{\bf Demonstration.} 
We start with the following essential proposition. It has been inspired
by  proposition 6.13 in \cite{ebin} which shows that the orbits of the action of \Dif\ on the manifold of Riemannian metrics are closed.
\begin{pro}\label{essent} 
Let $M$ be a connected manifold (not necessarily compact), let $\sui{h}$ be
a $C^k$-convergent sequence of pseudo-Riemannian metrics
on $M$ whose  limit  is $h_\infty$ $(k\geq 2)$, and let $\sui{\varphi}$ be a sequence of \diff s of
$M$ such that $(\varphi_n ^*h_n)_{n\in \mathbb N}$ is  $C^k$-conver\-gent
and tends to $g_\infty$. 
Let us  suppose there exist a convergent sequence $\sui{x}$ of points of
$M$ such that the sequence  $\suit {\varphi_n(x_n)}$ is  convergent,
and a real number $C>0$ such that $\forall \nN$,
$|D\varphi_n(x_n)|<C$.
Then the sequence  $\sui{\varphi}$ admits a   $C^k$-convergent
subsequence.
\end{pro} 
The notation $|.|$ refers to any fixed Riemannian metric. We need to
suppose that the sequence $\suit {\varphi_n(x_n)}$ is convergent to avoid situations like
translations in vector spaces where every point is sent to infinity
with bounded derivative. This is a typically non compact phenomenon. Hence, 
if we suppose $M$ compact we can simplify the statement.
\hop
{\bf Proof.} Denote by $x$ the limit of $\sui{x}$ and by $y$ the
limit of $\suit {\varphi_n (x_n)}$. Then the hypothesis 
 $|D\varphi_n(x_n)|<C$ induces the existence of a convergent sequence of frames $\sui{r}$ above
  $\sui{x}$, of limit  $r_x$, such that, restricting to a 
subsequence if necessary,
$$r'_n=D\varphi_n.r_n\ninf r'_y.$$
The frames   $r_x$ and  $r'_y$ define a linear map  $L$ from
$T_xM$ to $T_yM$ such that  $D\varphi_n(x_n)$ tends to $L$.\\
Next we want to define a map 
$\psi : M \rightarrow M$ 
such that a sub-sequence of $\sui{\varphi}$ converges to  $\psi$. 
Of course, we put  $\psi(x)=y$. Thanks to the following lemma, we are
going to extend  $\psi$ to a neighborhood of  $x$. We denote by ${\e}_n$
(resp. ${\e}'_n$) the exponential map of $h_n$ 
(resp. of $\varphi_n ^*h_n$)
and by ${\e}$  (resp. ${\e}'$) that of  $h_\infty$ (resp. $g_\infty$).
\begin{lem}\label{expo}
Let  $\sui{h}$ be a $C^k$-convergent sequence of pseudo-Riemannian metrics
whose limit is $h_\infty$. Then there exists an open  neighborhood
$U$ of the zero section of   $TM$  on which the exponential maps of
the metrics $h_n$ are all defined. Furthermore the sequence  $\sui{\e}$
converges  $C^{k-1}$ to $\e$ on every compact subset of $U$.
\end{lem}
{\bf Proof :} The $C^{k-1}$ convergence comes from the fact that the
exponential maps are solutions of differential equations whose 
coefficients, which  are the Christoffel symbols,  converge $C^{k-1}$.
We still have to check that the exponential maps are all defined at each
point on non trivial subsets of the tangent spaces.
To prove this, we adapt proposition 2.1 of 
A. Romero and  M. Sanchez in \cite{rom-san}.
Let  $v$ be a tangent  vector to $M$, $\sui{a}$ and  $\sui{b}$ be two
sequences of real numbers such that $[a_n,b_n]$ is the domain of $\gamma_n$ the
$h_n$-geodesic stemming from  $v$, and $a$ and $b$ such that
$[a,b]$ is that of $\gamma$ the $h_\infty$-geodesic. It is sufficient to show that 
$$\limsup a_n\leq a<0<b\leq \liminf b_n.$$
Put $A=\{t\in[0,b[| t<b_n$ for all $n$ but a finite number and $\gamma'_n(t)$ converge to $\gamma'(t)\}$.  Using that  solutions of
differential equations depend continuously on  the initial conditions \emph{
and} on  the coefficients of the equations, we see that $A$ is open and closed in $[0,b[$ and therefore $A=[0,b[$.  
$\Box$\hop
Let us  continue the proof of proposition \ref{essent}. The following relation holds:
\begin{rela}\label{permut} 
$\varphi_n\circ \e'_n=\e_n\circ D\varphi_n.$ 
\end{rela}
According to lemma  \ref{expo}, there is an open subset $U$ of $TM$ on
which  
$\e_n(\varphi_n(x_n))$ tends to  $\e(y)$ and  $\e'_n(x_n)$  tends to
$\e'(x)$. Restricting to an open subset if necessary, we can assert that  $\left(\e'_n(x_n) 
\right)^{-1}$ tends to   $\left(\e'(x)\right)^{-1}$. As  $x_n$  tends
to $x$,  there exists a neighborhood $V$ of
$x$ such that, for $n$ sufficiently large, 
$\left(\e'_n(x_n)\right)^{-1}$ is well defined on it.\\
Taking  $x'\in V$, one finds
$$\begin{array}{ll}\varphi_n(x')&=\varphi_n\left(\e'_n(x_n)\circ\left(\e'_n(x_n)\right)^{-1}(x')\right)\\ 
               &=\e_n\left(D\varphi_n(x_n)\circ\left(\e'_n(x_n) 
\right)^{-1}(x')\right), 
\end{array}$$
using relation \ref{permut}. However, the sequence of linear maps
$\suit{D\varphi_n(x_n)}$ converges to $L$ and therefore
$$\lim_{n\rightarrow\infty}\varphi_n(x')=\e(y)\circ L \circ (\e'(x))^{-1}(x'),$$ 
which we set to be $\psi(x')$. By construction, we see that the sequence $\sui{\varphi}$
converges,  $C^{k-1}$, to $\psi$ and that  $\psi$ is a local  $C^{\infty}$
diffeo\-morphism.
 Using the fact that  $D\varphi_n(x')$ tends to $D\psi(x')$, we repeat the same construction
to extend along the geodesics to the whole manifold. 
We have constructed a local $C^{\infty}$ diffeomorphism $\psi$  on  $M$.
The sequence $\sui{\varphi}$ converges, a priori, only $C^{k-1}$ to $\psi$. \\
Next we look at
the sequence of the connections, which is $C^{k-1}$ convergent.
We reproduce the computations done by D. Ebin  (\cite{ebin}, proposition 6.13).
In local coordinates we represent $\varphi_n$ by $\varphi^i_n(x^1,\dots,x^m)$.
Let ${}_n\Gamma^l_{ij}$ (resp. ${}_n\Gamma'{}^l_{ij}$) denote the Christoffel symbols of the metric $h_n$ (resp. of  $\varphi_n^*h_n$). We have
$${}_n\Gamma'{}^l_{ij}=\left(\frac{\partial \varphi_n^i}{\partial x^r}\right)^{-1}
                     \left(\frac{\partial \varphi_n^j}{\partial x^s}\right)^{-1}
                     {}_n\Gamma^t_{rs}
                     \left(\frac{\partial \varphi_n^l}{\partial x^t}\right)
                     -
                     \left(\frac{\partial^2 \varphi_n^l}{\partial x^s\partial x^r}\right)
                     \left(\frac{\partial \varphi_n^i}{\partial x^s}\right)^{-1}
                     \left(\frac{\partial \varphi_n^j}{\partial x^r}\right)^{-1},
$$
where ``$-1$'' indicates matrix inverse.
The convergence $C^{k-1}$ of the sequences of Christoffel symbols and the convergence $C^{k-2}$ of the first order derivatives of $\varphi_n$ implie the convergence $C^{k-2}$ of the second order derivatives, hence the convergence $C^{k}$ of the sequence of diffeomorphisms $\sui{\varphi}$.
\\
We still have to show that $\psi$  is a diffeomorphism, for this we
use the symmetry of the problem. Using that, asymptotically, the
$\varphi_n$ permute the induced volume forms of the metrics $h_\infty$
and  $g_\infty$, we see that the sequence $\sui{\varphi^{-1}}$  satisfies
also the hypothesis of the proposition and so we obtain the existence
of $\psi^{-1}$. This conclude the proof of proposition \ref{essent}
$\Box$ 
\vskip .25cm
\noindent{\bf Remarks.}\begin{enumerate}
\item
Proposition \ref{essent} can be used to prove that  on the
space of Riemannian metrics, the action is $C^2$-proper: indeed when $M$ is
compact and if the metrics are Riemannian, the hypotheses of the proposition are
always satisfied. 
\item
This proposition has an interesting interpretation: it
asserts that the sequences of diffeomorphisms we are interested in are
nowhere or everywhere divergent.\end{enumerate}
As in the proof of A. Zeghib, we first consider  the problem in the
linear case. Let us start the proof of theorem \ref{zorro} with the following evident lemma whose proof is left to the reader.
\begin{lem}\label{evid} 
If  $\sui{H}$ is a convergent sequence of quadratic forms of signature
$(m-1,1)$ on  $\mathbb R^m$ whose  limit is  $H$, then there exists a sequence 
of linear  maps   $\sui{M}$ of limit \emph{Id} such that $H_n={}^tM_n\,H\,M_n$.
\end{lem} 
This lemma shows that the linear approximately isometric systems
will have the same behavior as the non-equicontinuous (non bounded) sequences of
linear isometries. As a consequence, corollary  4.3 from \cite{zeghib} is
immediately transposable, and thus we have the following.
\begin{pro}\label{line}
Let  $(H_n,G_n,M_n)$ be a \emph{linear}  approximately isometric system
on $\mathbb R^m$. Then there exists a subsequence $\sui{N}$ such that
${AS}(0,N_n)$ is a lightlike hyperplane and moreover $\st {0,(N_n)}$  is the
set of strongly approximately stable vectors. 
Moreover, there exist
$C>0$ and a sequence of hyperplanes $\sui{P}$ such that we have the
following uniformity condition: $\forall
\nN,\ |{N_n}_{|_{P_n}}|<C$.
\end{pro}
Thanks to a trivialization of the tangent (measurable but continuous at
the neighborhood of the point considered) we associate to each
diffeomorphism a matrix field, $C_n(x)$.
We have 
$${PAS}\big(x,(\varphi_n)\big)={PAS}\big(0,(C_n(x))\big)={AS}\big(0,(C_n(x))\big).$$
Proposition  \ref{essent} tells us that if the sequence
$\sui{\varphi}$ is non-equicontinuous then, for all $x$, the sequence
$\suit{C_n(x)}$ is not bounded. By a diagonal process, we prove that
there exists a subsequence $\sui{\phi}$ of $\sui{\varphi}$ 
and a dense subset $M'$ of $M$ such that, for all $x\in M'$, $PAS(x,\phi_n)$ is a hyperplane of
$T_xM$ that we will denote by $P_x$. 
Keeping the same notations for
the exponential maps as before, we can prove the following.
\begin{pro} For all  $x\in M'$, there exists a  neighborhood  $V_x$ of $0$
in $T_xM$ such that  the hypersurface
$\mca I_x=\e'(P_x\cap V_x)$  is  geodesic and such that for all  $y\in \mca I_x$,
$T_y\mca I_x=AS(y,\phi_n)$.
\end{pro}
{\bf Proof.}
 Let  $\sui{P}$ be the sequence of hyperplanes in $T_xM$ given by
proposition \ref{line} applied to the  $C_n(x)$.
We are going to use the fact (cf. lemma
\ref{expo}) that the sequences of exponential maps converge. Let $y$
 be a point of $\mca I_x$, we first
 write the differential of the $\phi_n$ at this point using the
 exponential maps of the $g_n$ that we denote, as before, by $\e'_n$.
\begin{eqnarray*}
D_y\phi_n&=D_{\left(\e'_n(x)\right)^{-1}(y)}\Big(\phi_n \circ
                \e'n(x)  \Big)
                \circ    D_y  \left(\e'_n(x)\right)^{-1}    \\ 
                &=D_x\Big( \e_n(\phi_n(x)) \circ \phi_n\Big) \circ
                D_y\left(\e'_n(x)\right)^{-1},
\end{eqnarray*}
according to relation \ref{permut}.
We set  $P'_n=T_y\left(\e'_n(P_n\cap V_n)\right)$, with $V_n$ a bounded
neighborhood
of  $0$ in $T_xM$ inside the domain in which  $\e'_n$ is defined.
It is clear that  $P'_n\rightarrow T_y\mca I_x$, but also that  
$D_y\left(\e'_n(x)\right)^{-1}. P'_n=P_n$. Consequently, as the
sequence $\sui{\e'}$ is convergent, $|{D\phi_n}_{|P'_n}|$ is bounded.
Hence, $T_y\mca I_x\subset AS(y,\phi_n)$. However, thanks to
proposition \ref{essent} (following the proof of the fact 6.4 of \cite{zeghib}), we can prove 
that $AS(y,\phi_n)$ is at most
a hyperplane therefore we have equality.

Now, we have to prove that this surface is geodesic.
Let ${c}$ be a path on $\mca I_x$ between $x$ and $y$.
The neighborhood $V_n$ being bounded, we see that
the sequence $\sui{d}$, defined by
$d_n=\phi_n\circ c$, is with $C^1$ bounded variation.
We denote by $\tau_n$ the parallel transport along $c$ for
the  metric $h_n$ and by $\tau'_n$ the parallel transport
along  $d_n$ for the metric $g_n$.
Finally, we denote by $P''_n$ the image of $P_n$ by $\tau_n$ and by
$P''$ the limit of this sequence.
We have:
$$D_{y}\phi_n=\tau'_n\circ  D_x\phi_n\circ\tau^{-1}_n.$$
Knowing that the sequence of paths  has bounded  variation and the
sequences of metrics are $C^1$ convergent, it becomes clear that we
have a new sequence of stable hyperplanes that is $P''$ is
approximately stable. Therefore $P''=AS(y,\phi_n)$ and so
$P''=T_y\mca I_x$. The  surface is therefore geodesic.
$\Box$\hop
The size of the neighborhood $V_x$ depends  essentially on the
 domains on which the exponential maps are defined.
Moreover, two plaques
 $\mca I_x$ and  $\mca I_y$ are either disjoint or tangent.
Thus, the hypotheses of the following lemma  (cf. \cite{zeghib} and
\cite{ze2}) are fulfilled. 
\begin{lem}
Let $M$ be a compact manifold endowed with a torsion free connection
and an auxiliary norm $|.|$ on $TM$. Let $M'$ be a dense subset of $M$ and
suppose given a real number $r$ and for $x\in M'$ a hyperplane $P_x\subset
T_xM$ and let  $\mca{I}_{x,r}=exp(P_x\cap B_x(r))$, where $B_x(r)$ is
the ball of $T_xM$ centered at $0$ and with radius $r$. Also, suppose
that $\mca{I}_{x,r}$ is geodesic  and that if two plaques
$\mca{I}_{x,r}$ and  $\mca{I}_{y,r}$ intersect at some point, then
they are tangent at that point (and hence by geodesibility, the
intersection $\mca{I}_{x,r}\cap\mca{I}_{y,r}$ is open in both
$\mca{I}_{x,r}$ and  $\mca{I}_{y,r}$). Then the geodesic plaques
$\mca{I}_{x,r}$ extend to a Lipschitz geodesic foliation of $M$.
\end{lem}
Now, to  finish the proof of the first part, we have to show that
the tangent to this foliation is still approximately stable.
Moreover, we can prove that
$AS(\phi_n)=PAS(\phi_n)$, and also  that  the set of strongly approximately stable vectors
is the tangent bundle to the $1$ dimensional geodesic lightlike
foliation $\st{\phi_n}$. The proof of those three points following exactly Zeghib's proof, we refer to \cite{zeghib}.

We only give the main idea of the proof of the second part of the theorem. It is exactly the
 same as the original case (proposition 9.5 of \cite{zeghib}). It uses the fact that
 if $v\notin AS(\phi_n)$ then there exists a sequence of real numbers $\sui{\alpha}$ and a sequence of vectors
$\sui{u}$ such that $D\phi_n v=\alpha_n u_n$. Of course, $\alpha_n\rightarrow\infty$ and therefore 
$|D\phi_n^{-1}(u_n)|\rightarrow 0$. The limit of the sequence of vectors is consequently in 
$AS^\perp(\phi_n^{-1})$. 

Let us show the last part of the theorem.
Let $x$ be a point of $M$,  and let $\sui{v}$ be a convergent sequence of
vectors which  tends to $v$ and such that $D\phi_n(v_n)\rightarrow 0$,
i.e. $v\in \st{\phi_n}$. Let us consider
$\Delta_n=\sigma'_n(\e_n(x).v_n)-\sigma'_n(x)$, we have 
$$\begin{matrix} 
\Delta_n&=\sigma_n(\phi_n(\e_n(x).v_n))-\sigma_n(\phi_n(x))\\ 
        &=\sigma_n(\e'_n(D\phi_n\, v_n))-\sigma_n(\phi_n(x)), 
\end{matrix}$$
using relation \ref{permut}. If we take  $\sui{v}$ such that 
$v_n\rightarrow v$ and $D\phi_n(v_n)\rightarrow 0$, we have
$$\lim_{n\rightarrow\infty}\e'_n(D\phi_n(v_n))=\lim_{n\rightarrow\infty}\phi_n(x),$$ 
therefore
$$\sigma_n(\e'_n(D\phi_n\, v_n))-\sigma_n(\phi_n(x))\ninf 0$$ 
and therefore $\Delta_n$ also converges to $0$. However, $\sigma'_n(\e_n(x).v_n)$ converges
to $\sigma'(\e(x).v)$ and $\sigma'_n(x)$ converges to $\sigma'(x)$.
Therefore, we can assert that   $\sigma'$ is constant along the
$g_\infty$-geodesic stemming from $v$, i.e. $\st{\phi_n}$ according to
what preceded.
$\Box$ \hop
In \cite{DGr}, G. D'Ambra and M. Gromov conjectured that if $M$ is compact and simply connected then the action is $C^2$-proper. Theorem \ref{zorro} solves the $3$-dimensional case.
\begin{coro}\label{coro} 
Let $M$ be a $3$ dimensional compact manifold not covered by $\mathbb
R^3.$ Then the action of the group of diffeomorphisms on $\mca L(M)$ is proper 
for the $C^2$ topology. 
\end{coro}
{\bf Proof :} According to  \cite{ze2}, these manifolds do not possess 
  lightlike geodesic foliations of  codimension $1$ therefore theorem \ref{zorro}
 finishes the proof. $\Box$ \\
This statement can be found  without proof in \cite{zeghib}.
\vskip .2cm\noindent
Since any Lorentzian manifold possess a negative line-field, we see that the metrics $h_\infty$ and $g_\infty$ of theorem \ref{zorro} have two nowhere collinear line-fields. It trivially gives a topological obstruction. 
For example, thanks to  the famous theorem of J.F. Adams (see \cite{jfa}) about vector fields on spheres, we have:
\begin{coro}
If $n\equiv 2$ mod $4$, then the action of $\mathrm{Diff}(S^{n-1})$ on $\mca L(S^{n-1})$ is proper.
\end{coro}
\section{On the  Torus.} 
\subsection{Localization of the non properness.}\label{P}
There exist (flat) metrics on the torus with non compact isometry group, consequently the action  of the group of \diff s of the torus on $\mca{L}(\mathbf T^2)$ is not proper. A priori, the  approximately isometric systems $(h_n,g_n,\varphi_n)$ can have a very different shape. The point is to understand it.
Even if every Lorentzian surface admits two  geodesic
lightlike foliations of (co)dimension $1$, theorem \ref{zorro} will enable us
to localize precisely the non properness of the action
 by showing that the limit metrics 
$h_\infty$ and  $g_\infty$ are flat. We first show the following result.
\begin{theo}\label{bonne}
A Lorentzian metric on the torus $\mathbf T^2$ whose curvature is
constant along one of its lightlike foliations is flat.
\end{theo}
Let us remark that, of course, the lightlike condition is necessary as
 the case of the Clifton-Pohl torus shows (cf. \cite{carr-roz}). There also exist
local counterexamples: if we take the metric $h(x,y)=dxdy+2\,xy^2\,dx^2$, its
curvature is $K(x,y)=-x$ which is constant along the lightlike direction $\partial_y$.\hop
{\bf Proof:} We call $\ff$  the foliation along which the curvature is
constant. According to the Gauss-Bonnet theorem \cite{nomizu}, the
curvature has to vanish somewhere.  Therefore there exists a leaf $F$
of $\ff$ such that the curvature, $K$, vanishes along $F$. Two cases have
to be considered: the leaf can be compact or not.
Let us suppose the leaf is not compact. If the foliation is conjugated
to a linear foliation then $F$ is everywhere dense and 
$h_\infty$ is clearly flat. If this is not the case, $\ff$ has some
compact leaves (cf. \cite{HH}). We can thus reduce the problem to a
foliation on an annulus. Therefore we can state (cf. theorem 4.2.15 \cite{HH}) that the leaf $F$ is going
to accumulate on two compact leaves. We denote by $A$ the Lorentzian
annulus delimited by those leaves.
All the leaves of $\ff$ contained in $A$ are going to accumulate on the
boundary of $A$. Therefore this annulus is flat.
\begin{lem}
If $A$ is a flat Lorentzian annulus with lightlike boundary and
$\ff$ is the lightlike foliation tangent to the boundary, then all the
leaves of $\ff$ are compact.
\end{lem}
{\bf Proof:} Let $\widetilde A$ be the universal cover of $A$
( topologically $\widetilde A$ is a stripe) and $\gamma$ be a  generator of 
$\pi _1(A)$.
The flat Lorentzian manifold $\widetilde A$ is developed by a local diffeomorphism
 $D:\widetilde A\longrightarrow{\mathbb R}^2$ (where ${\mathbb R}^2$
is endowed with the Lorentzian metric $xy$) which satisfies the
following equivariant condition : 
$$(*)\quad \quad D\circ\gamma=\gamma'\circ D,$$
 where $\gamma'$  is a Lorentzian isometry 
(the holonomy of $\gamma$). 
The foliation $\ff$ is lifted to a foliation 
$\widetilde{\ff}$ of $\widetilde A$. We can suppose that $\widetilde{\ff}$
is the pull-back  by $D$ of the foliation  of ${\mathbb R}^2$ by  
horizontal lines.
 
Translating the situation if necessary, we have two cases~:\\
a) $\gamma'=
\begin{pmatrix}
\lambda &0\\ 
0&1/\lambda
\end{pmatrix},$ with $\lambda \in \mathbb R^*$\\
b) $\gamma'$ is an horizontal translation.

The case a) is impossible; in fact, the continuous function $f(x,y)=xy$
is $\gamma'$-invariant and, on the closed superior and inferior
half-planes, its only possible local extremum is $0$.
According to the condition $(*)$ and the fact that $D$ is a local diffeomorphism,
the continuous  function  $g=f\circ D$ is $\gamma$-invariant
and its only local  extremum is $0$. Therefore,  it descends to $A$
as a  constant  function everywhere equal to $0$, which is not possible.\\
Therefore the only remaining case is  b). In this case, the annulus
$A$ (for the same 
 affine structure) admits a  flat Riemannian metric, which implies
that $D$ is a diffeomorphism (the structure is complete as
 $A$ is compact). Consequently, $\widetilde A$ can be directly seen as
a horizontal stripe of
 ${\mathbb R}^2$ and we obtain $A$ by quotienting by a
horizontal translation. Hence, all the leaves of  $\ff$
 are closed.
$\Box$ \hop
We deduce from this that if  ${\cal F}$ has no dense leaf (we can
assume this) then a leaf $F$,  along which  the curvature is zero, has to be 
 compact.
Moreover, let us show that $F$ is geodesically complete.
Let  $\varphi $ be a germ of  diffeomorphisms at  $0\in{\mathbb R}$
which  
generates the holonomy of $F$. According to the article of
Y. Carri\`ere and L. Rozoy [C-R], the completeness of
$F$ is  characterized by the divergence of the two
series $\sum _0^{\pm\infty}\varphi'(0)^k$
i.e. by $|\varphi'(0)|=1$. But, if $|\varphi'(0)|\ne 1$ then
$F$ is attractive (or repulsive according to the choice of an
orientation).
A leaf attracted by $F'$ would be non compact and the curvature would
vanish along it. That contradicts what we just proved.
Hence every leaf of $\ff$ along which the curvature vanishes is a
 compact and complete geodesic.\\
If the curvature is not constant around $F$ we can consider, using
 a transverse curve, the closest leaf of $\ff$ with zero curvature. Again we
 obtain a Lorentzian annulus $A'$ (we can also take the closure
 of a connected component of  $\mathbf T^2\setminus\{K^{-1}(0)\}$).
The connected components of the boundary of $A'$ are complete lightlike geodesics and the
 curvature does not vanish inside $A'$.
We are going to use the following version of the Gauss-Bonnet theorem.
\begin{lem}\label{GB}
Let $A$ be a Lorentzian annulus with lightlike boundary
$\Gamma=\gamma_1\cup\gamma_2$.
Let $\gamma_i :[0,1]\rightarrow \Gamma $, $i\in \{1,2\}$ be a
parametrization of the boundary. Let $Z_i$ be a vector field along
$\gamma_i$, tangent to $\gamma_i$ and parallel. We define $\lambda_i$ as
the proportionality coefficient between  $Z_i(0)$ and  $Z_i(1)$. Then
we have\,:
$$\int_A K dv_h=\ln(\lambda_1/\lambda_2).$$
If the boundary is made of  complete geodesics,  we  have:
$$\int_A K dv_h=0.$$
\end{lem}
{\bf Proof:} Let  $X$ be a nowhere vanishing  lightlike  vector field
tangent to the boundary,  $X^0$ be another lightlike vector field such
that  $h(X,X^0)=1$ and  $\omega$ be the 1-form defined by   
$\omega(v)=h(\nabla_vX,X^0)$, where $\nabla$  is the Levi-Civita
connection of the metric.
Then:
$$\begin{array}{lll} 
d\omega(v,u)&=&v.\omega(u)-u.\omega(v)-\omega([v,u])\\ 
            &=&v.h(\nabla_uX,X^0)-u.h(\nabla_vX,X^0)-h(\nabla_{[v,u]}X,X^0)\\ 
            &=&h(\nabla_v\nabla_uX,X^0)+h(\nabla_uX,\nabla_vX^0)-h(\nabla_u\nabla_vX,X^0)\\ 
            & &-h(\nabla_vX,\nabla_uX^0)-h(\nabla_{[v,u]}X,X^0)\\ 
            &=&h(R(v,u)X,X^0)+h(\nabla_uX,\nabla_vX^0)-h(\nabla_vX,\nabla_uX^0). 
\end{array}$$
Contrarily to the Riemannian case and to the case where $X$ is
timelike (see \cite{nomizu}), we do not have automatically
$h(\nabla_uX,\nabla_vX^0)=0$.
Nevertheless, we still have  $h(\nabla_uX,X^0)=0$ and so
$h(\nabla_uX,\nabla_vX^0)=h(\nabla_uX,X^0)h(X,\nabla_vX^0)$. 
Derivating the equality $h(X,X^0)=1$, we obtain immediately
$h(\nabla_uX,X^0)=-h(X,\nabla_uX^0)$.
Therefore we have
$$\begin{array}{ll} 
h(\nabla_uX,\nabla_vX^0)&=h(\nabla_uX,X^0)h(X,\nabla_vX^0)\\ 
&=-h(X,\nabla_uX^0)(-h(\nabla_vX,X^0))\\ 
&=h(\nabla_vX,\nabla_uX^0). 
\end{array}$$
Hence $d\omega(X,X^0)=h(R(X,X^0)X,X^0)=K$ and so  $d\omega=Kdv_h$.
According to  Stokes' theorem, we have\,:
$$\int_A Kdv_h=\int_\Gamma \omega,$$
where  $\Gamma=\gamma_1\cup\gamma_2$ is the boundary of $A$.
Let $Z_1$ be a  parallel vector field tangent to  $\gamma_1$. Then 
$Z_1=\mu X$, and we obtain 
 $$\nabla_{Z_1}Z_1=(Z_1.\mu)X+\mu\ \omega(Z_1) X=0.$$
Thus we have $\omega(Z_1)=(-Z_1.\mu)/\mu$, and therefore
$\int_{\gamma_1}\omega=-\ln\mu(0)+\ln\mu(1)=-\ln \lambda_1$.
Doing the same with $\gamma_2$ yields
$$\int_{\Gamma}\omega=\ln(\lambda_2/\lambda_1).$$
$\Box$\hop
We can now finish the proof of the theorem.
On the annulus $A'$ we had previously, we apply lemma \ref{GB} and find
$\int_{A'} K=0$. However 
$K$ does not vanish inside  $A'$ and thus we have a contradiction and
therefore $K=0$  on all $\mathbf T^2.$
$\Box$ \\
This result applies directly to our problem. If, as in the introduction, we denote by $\mca F(\mathbf T^2)$ 
the set of flat metrics of the torus, the following is true.
\begin{coro}\label{plat}
The action of {\rm Diff}$(\mathbf T^2)$ on $\mca L(\mathbf T^2)-\mca F(\mathbf T^2)$ is $C^k$-proper, 
$k\geq 2$. 
\end{coro}
{\bf Proof:} It is clear that the statement is a consequence of the following stronger proposition:
if $(h_n,g_n,\varphi_n)$ is a  $C^2$ approximately isometric system on
$\mathbf T^2$, then the metrics  $h_{\infty}$ and
$g_{\infty}$ are flat.
 According to theorem \ref{zorro}, $\st{\varphi_n}$ and  $\ist{\varphi_n}$ are
lightlike Lipschitz line-fields (actually, in dimension $2$ they are
automatically smooth). We can apply the second part of the theorem to
the curvature functions. Indeed, if we denote by  $K_n$ the curvature
of $h_n$ and  $K'_n$ the one of $g_n$, these functions satisfy the
relation 
$K_n\circ\varphi_n=K'_n$ and the $C^2$-convergence of the sequences of
metrics  implies the convergence of the sequences of curvature functions.
Hence, the curvatures of both $h_\infty$ and  $g_\infty$ are constant
along a lightlike foliation. According to the previous theorem
these metrics are flat.$\Box$\hop
This corollary will be of  precious help to study the manifestations
of the non properness of the action of \Dif\ on $\mca L(\mathbf T^2)$.
But, on the Klein bottle, it entails the  following interesting result. 
\begin{coro}\label{bouteille}
Let  $K$ be the  Klein bottle. The action of  {\rm Diff}$(K)$ on  $\mca{L}(K)$
is  proper.
\end{coro}
{\bf Proof:} 
We are going  to show that the sequences of diffeomorphisms  that compose the
approximately isometric systems of the torus can never descend to the
Klein bottle. For this we first show the
\begin{lem}\label{identite}
Let $(h_n,g_n,\varphi_n)$ be a $C^2$-approximately isometric system
on $\mathbf T^2$ and let  $f_n$ be the lift  to  $\mathbb R^2$ of  $\varphi_n$.
Then there exists a sequence $\sui{N}$ of affine  maps of $\mathbb R^2$ such that the sequence 
$\sui{N_n\circ f}$ converges to the identity map.
\end{lem} 
{\bf Proof:} As we permit right or left composition by \diff s, we can suppose that the lifts 
of  $h_{\infty}$ and of
$g_{\infty}$, that we denote by $H$ and $G$, are quadratic
forms and even that $H=G$. 
 We denote by $H_n$ (resp. $G_n$)
the lift of the metric $h_n$ (resp. $g_n$).
Of course, $f_n^*H_n=G_n$  and converges to $G$. We denote by $x_n$ a point of $\mathbb R^2$ that realizes the 
maximum of the norm of the second derivative of $H_n$ (seen as a map from $\mathbb R^2$ to $GL(2,\mathbb R^2)$).
Let $M_n$
be the   linear map  $Df_n(f_n^{-1}(x_n))$.
We know that ${}^t\!M_nH_n(x_n)M_n$ converges to $H$, hence from lemma
 \ref{evid}, $M_n$ is close to an isometry of $H$. Let us choose a basis of $\mathbb R^2$ such that the metrics are
given by the matrices
$$H=\begin{pmatrix}
0&1\\
1&0
\end{pmatrix}\qquad \mrm{and} \qquad
H_n(x)=
\begin{pmatrix}
a_n(x)&b_n(x)\\
b_n(x)&c_n(x)
\end{pmatrix}.
$$
We can suppose that $M_n$ is close to 
$$\begin{pmatrix}
1/\lambda_n&0\\
0&\lambda_n
\end{pmatrix},$$
with $\lambda_n\rightarrow \infty$. \\
As $H_n\rightarrow H$, for the $C^2$ topology, it is easy 
to see that  
$D^2\ ({}^t\!M_n H_n(x) M_n)\rightarrow 0$ 
uniformly if and only if 
$D^2(\lambda^2_nc_n)\rightarrow 0$ uniformly. As $f_n^*H_n\rightarrow H$, we have 
$D^2(\lambda^2_nc_n(x_n))\rightarrow 0$ and as $x_n$ realizes the maximum of the second derivative we have 
the uniform convergence to $0$.
Now we use again that   $f_n^*H_n\rightarrow H$ to show that 
$D\,({}^t\!M_n H_n(x_n) M_n)\rightarrow 0$ and that 
${}^t\!M_nH_n(x_n)M_n\rightarrow H$. We deduce from this that $G'_n=M_n^*H_n\rightarrow H$ 
for the $C^2$ topology. We can define a sequence of translations $\sui{\tau}$, such that 
$\tau_n (M_n^{-1}(x_n))=0$. The sequence $G''_n=\tau_n^*G'_n$ still converges. We denote by $N_n$ the 
map $(M_n\circ \tau_n)^{-1}$.
Therefore
$(N_n\circ f_n)^*G''_n\rightarrow G$. However, the
sequence of the derivatives of $(N_n\circ f_n)$ at $f_n^{-1}(x_n)$ is bounded and 
$\left(N_n\circ f_n\right)(f_n^{-1}(x_n))=0$.
We recognize the hypothesis of proposition \ref{essent}. Extracting a
subsequence if necessary, we can say that the sequence
$\suit{N_n\circ f_n}$ is $C^2$-convergent. Its limit is an
isometry of $H$ tangent to a translation and therefore it is a translation. It is not hard, now,
to obtain the desired sequence.
 $\Box$\hop
We continue the proof of the corollary. Taking $n$
sufficiently big, we write $f_n=M_n\circ \varepsilon_n$ 
with $M_n$ linear and 
 $\varepsilon_n$ a \diff\ close to the identity. Let $\sigma$ be the
``antipodal'' map  defined by  $\sigma(x,y)=(x+1/2,-y)$. The fact that
$f_n$ actually descends to the Klein bottle can be written
$f_n\circ\sigma=\sigma\circ f_n$.
We have 
\begin{eqnarray*}
M_n\circ\varepsilon_n\circ\sigma=&\sigma\circ M_n\circ\varepsilon_n\\
M_n\circ\varepsilon_n\circ\sigma\circ\varepsilon_n^{-1}=&\sigma\circ M_n,
\end{eqnarray*}
these equalities being on the torus. We set 
$\sigma_n'=\varepsilon_n\circ\sigma\circ\varepsilon_n^{-1}$.
Of course, $\sigma$ and $\sigma_n'$ are close.
Then we compute  $(M_n\circ\sigma_n'-\sigma\circ M_n)(x,0)$ and
$(M_n\circ\sigma_n'-\sigma\circ M_n)(0,y)$ which have to be in $\mathbb
Z^2$. It gives immediately that $M_n$ has to be diagonal. Using again
the fact that  $\sui{f}$ descends to the quotient and is not
equicontinuous, we see that it is impossible. The sequence cannot descend to the Klein bottle.
$\Box$
\subsection{Action of {\rm Diff}$(\mathbf T^2)$  on $\mca F(\mathbf T^2)$.}\label{F}
As we announced in the introduction, we are going to 
prove that the action of the group of diffeomorphisms on
 the set of flat metrics of volume $1$, denoted by $\mca F_1(\mathbf T^2)$, is ergodic 
(or more correctly that it has the same dynamical properties as  a given ergodic action).
If we want to talk about ergodicity, we first have to reduce the problem to
 a finite dimensional one. This reduction is done by the following proposition. 
\begin{pro}
The action of \emph{Diff}$_0(\mathbf T^2)$ on $\mca F_1(\mathbf T^2)$ is proper and 
$\mca F_1(\mathbf T^2)/\mrm{Diff}_0(\mathbf T^2)$ is diffeomorphic to $SL(2,\mathbb R)/SO(1,1)$.
\end{pro}
{\bf Proof:} Let us consider an approximately isometric system. To see the properness, we  
improve  lemma \ref{identite} by noting that the linear part of the affine maps involved
 can be chosen in $GL(2,\mathbb Z)$. As Diff$(\mathbf T^2)/$Diff$_0(\mathbf T^2)$ is isomorphic to 
$GL(2,\mathbb Z)$, the diffeomorphisms of an approximately isometric  system cannot stay in Diff$_0(\mathbf T^2)$.\\
Let us show now that the quotient is diffeomorphic to $SL(2,\mathbb R)/SO(1,1)$, the space of quadratic forms of signature $(1,1)$ and of determinant $1$. 
Let $h$ be a flat metric. As $h$ is complete, there exists a \diff \ $\psi$ such that $\psi^*h$ is a quadratic form.
As Diff$(\mathbf T^2)/$Diff$_0(\mathbf T^2)$ is isomorphic to 
$GL(2,\mathbb Z)$, there exists a unique element $N$ of $GL(2,\mathbb Z)$ such that $\psi$ is equivalent to $N$ modulo 
Diff$_0(\mathbf T^2)$. Thus, $\psi\circ N^{-1}$ lies in Diff$_0(\mathbf T^2)$ and  $(\psi\circ N^{-1})^*h$ is the only  quadratic form in the Diff$_0(\mathbf T^2)$-orbit of $h$.
$\Box$
\vskip .3cm\noindent
This proposition implies that the properties of the action of Diff$(\mathbf T^2)$ on 
$\mca F(\mathbf T^2)$ are the same as those of the left action of Diff$(\mathbf T^2)/$Diff$_0(\mathbf T^2)$ on the set of quadratic forms of signature $(1,1)$. We restrict ourself to the action of Diff$^+(\mathbf T^2)$. As Diff$^+(\mathbf T^2)/$Diff$_0(\mathbf T^2)$ is isomorphic to $SL(2,\mathbb Z)$,
 the point is now to understand  the action of $SL(2,\mathbb Z)$ to $SL(2,\mathbb R)/SO(1,1)$ and to show that it is ergodic.
\vskip .2cm \noindent
The proof will be complete by noting  that this action has the same dynamical properties as the right action of $SO(1,1)$ on  $SL(2,\mathbb Z)\backslash SL(2,\mathbb R)$. This is due to the fact that the lifts to $SL(2,\mathbb R)$ of the invariant sets and functions  of those actions are the same.
This third action is well known, it is the action of the geodesic flow of the unitary tangent bundle 
of the modular surface 
$SL(2,\mathbb Z)\backslash \mathbb H^2$, where $\mathbb H^2$ is the hyperbolic plane.
This manifold has finite volume and the action is known to be ergodic. This implies that our action is also ergodic.
 In particular, almost all flat orbits of $\mca F_1(\mathbf T^2)$ are dense in $\mca F_1(\mathbf T^2)$. 
Therefore we have proved:
\begin{theo}\label{ergod}
The action of \emph{Diff}$(\mathbf T^2)$ on  $\mca F_1(\mathbf T^2)$ is ergodic.
\end{theo}
It implies that, on the set of flat metrics, there are no non constant continuous function invariant by rescaling and by the action of Diff$(\mathbf T^2)$.
\\
Now, to complete the description, we may wonder which are the closed flat orbits. 
We can give the following proposition.
\begin{pro}\label{orbi}
The orbit of a flat metric for which  all lightlike curves are
closed is $C^2$-closed. Moreover an orbit is $C^2$-closed if and only if it is $C^\infty$-closed. 
\end{pro}
{\bf Proof:} It is clear that non-closed orbits always contain an approximately isometric system whose shape is $(h,g_n,\varphi_n)$, ie with $\sui{h}$ constant. We are going to see that such systems do not exist in the  case above. 
We first give the following lemma.
\begin{lem}\label{suite}
Let $H$ be a quadratic form of signature $(1,1)$ on $\mathbb R^2$ and
$\sui{M}$ be a non bounded sequence of automorphisms such  that $\sui
{{}^tM_n\,H\,M}$ is convergent. Then there exists a sequence
$\sui{I}$ in Is$(H)$  such that the sequence  $\sui{I_n^{-1}\,M}$ is
bounded. Moreover there exists an $H$-lightlike vector $v$ such that $M_n^{-1}v$ tends to $0$.
\end{lem}
{\bf Proof:} The first part of this lemma is trivially deduced from  lemma  \ref{evid},
and the second one is an immediate consequence of the shape of the elements of $O(1,1)$.
$\Box$ 
\hop
We suppose that such an approximately isometric system exists. Then by lifting the system  to $\mathbb R^2$, we can replace the diffeomorphisms by elements of $GL(2,\mathbb Z)$, thanks to lemma \ref{identite}.
Set  $v=(\tau,1)\in\istb{\varphi}$ and let $H$ be a Lorentzian quadratic form and $\sui{M}$ be a non bounded sequence in $GL(2,\mathbb Z)$ such that
$\sui{{}^tM_n\,H\,M}$ is convergent. We set
$$M_n^{-1}=\matg{q}{p}{r}{s}.$$  
According to lemma \ref{suite}, we have $M_n^{-1}=P_n^{-1}\,I_n^{-1}$,
where $I_n\in \mrm{Is}(H)$  and the sequence $\sui{P}$, of elements of
$GL(2,\mathbb Z)$, is convergent. The sequence 
 $\sui{I}$ not being bounded, there exists a subsequence of $\sui{M}$
such that $M_n^{-1}(v)\rightarrow 0$. We thus have
$q_n\,\tau+p_n\rightarrow 0$. If $\tau\in\mathbb Q$ this sequence is
stationary (cf. \cite{HW}) and so is the sequence $r_n\,\tau+s_n$
therefore the matrix cannot be invertible. Hence, $\tau\in\mathbb
R\setminus \mathbb Q$ and $h$ has no compact lightlike geodesics and thus we have a contradiction.

The last assertion is immediate because the non $C^2$-closed orbits involve quadratic forms and linear maps and therefore the convergence can be supposed $C^\infty$. 
$\Box$\hop
{\bf Remarks.}\begin{enumerate}
\item After this proof, we might wonder if, on the torus, the approximately stable foliations of approximately isometric system  always have dense leaves.
This turns out to be false as shown by the example of \cite{these} p. 48.
\item 
A priori the statement of proposition \ref{orbi} is not optimal. Actually, we can prove that for any flat metric $h$ with non 
compact lightlike curves there exists a  \emph{non equicontinuous} sequence of diffeomorphisms $\sui{\phi}$ such that the sequence $\suit{\phi_n^*h}$ is $C^\infty$-convergent (see \cite{these} theorem IV.3). 
Anyway, it is not enough to prove that its orbit is not closed, we still have to compare $h$ with the limit of $\suit{\phi_n^*h}$.  
\item\label{anosov}
Let us give an example of a non flat orbit which is not $C^1$-closed.
Let $A$ be  the Anosov map, that is the torus diffeomorphism induced
by the matrix
$\begin{pmatrix}
2&1\\
1&1
\end{pmatrix}$. 
Let $g_A$ be an $A$-invariant Lorentzian (flat) metric on $\mathbf T^2$, and 
let $f$ be a non constant function of the torus. Let $X$ (resp. $Y$) be a non zero vector field of
the contracting (resp. dilating) direction of $A$.
We disturb  $g_A$  along the lightlike
contracting direction $X$ (using $Y$) by defining:
$$h=g_A+f\,(Y^\flat \otimes Y^\flat),$$
The metric is still Lorentzian.
The sequence $\suit{A^{n\,*}h}$ is $C^1$-convergent and its limit is
$g_A$, but $h$ is not flat.
Hence, there exist orbits which are not $C^1$-closed but which are
$C^2$-closed.
\end{enumerate}
\section{The set of metrics without isometries.} \label{O}
If $M$ is a closed manifold, we recall that, in the Riemannian case, the set of metrics without non trivial isometries, that we will denote by $\mca L_T(M)$, is an  open dense subset of the space
of Riemannian metrics endowed with the $C^\infty$ topology, see
\cite{ebin}.
This result is a corollary of the slice theorem of Ebin which is now clearly not
true in the Lorentzian case. Anyway, adapting the proof  
of Ebin (see proposition IV.4.2 of \cite{these}) we still have:
\begin{pro}
The set $\mca L_T(M)$ of metrics without non trivial isometries is  dense in $\mca{L}(M)$.
\end{pro}
The main point is to show that this set is open. Unlike to the
Riemannian case the answer is different according to which manifold we
consider. 
\subsection{A positive result.}
Thanks to corollary \ref{plat} and 
the theorem of D. Ebin, we can give the following
result. 
\begin{theo}\label{ouvert} 
Let $M$ be the $2$-torus or a  manifold such that  {\rm Diff}$(M)$ acts
properly\footnote{see corollary \ref{coro} and corollary \ref{bouteille}}, for the
$C^\infty$-topology, on $\mca L(M)$. Then  $\mca L_T(M)$ is  an open (dense) set.
\end{theo}
{\bf Proof:} Let $M$ be a manifold as in the statement of the theorem. 
Let $\sui{h}$ be a $C^\infty$ convergent sequence of Lorentzian metrics on
$M$, such that for all $\nN$ the group of isometries of
$h_n$, Is$(h_n)$, is non trivial. We are going to show that $h_\infty$, the limit
of the sequence of metrics, also has a non trivial isometry.
For this, let $\varphi_n$ be a non trivial element of  {Is}$(h_n)$.
If $\sui{\varphi}$ is convergent and if its limit is not the identity
map, then we have the result.
If the action is proper the sequence has to be equicontinuous.
If we are on the torus and if the sequence is not equicontinuous then
corollary \ref{plat} tells us that $h_\infty$ is flat and therefore
has a big isometry group.
Consequently, in both cases, we can consider that
the sequences of isometries we consider converge. 
Furthermore, thanks to corollary \ref{plat}, we know that the metrics of the torus with non compact isometry groups are flat (more precisely they are invariant under an Anosov map like the metric $g_A$ of the remark \ref{anosov}, page \pageref{anosov}). Consequently, the isometry groups of the metrics $h_n$ always contain a compact subgroup. 
Therefore we can assume that the $\varphi_n$ are of finite order $p_n$.
Now, we want to prove that there exists a sequence of integers  $\sui{k}$
such that $\sui{\varphi^{k_n}}$ admits a convergent subsequence whose 
limit is different from the identity.\\
Let us suppose it is not true, then for any sequence $\sui{k}$ of integers,
$\varphi_n^{k_n}\rightarrow$~Id.
Let $\alpha$ be a Riemannian metric with trivial isometry group. Then we define
$$\alpha_n=1/p_n\sum_{k=1}^{p_n}\left(\varphi_n^k\right)^*\alpha,$$
where $p_n$ is the order of $\varphi_n$. We have, by construction,
$\varphi_n^*\alpha_n=\alpha_n$. 
We are going to show that $\alpha_n\rightarrow\alpha$, for the
$C^\infty$-topology.
For this denote by $\|.\|_l$ the norm which gives the $C^l$ topology. We
have:
$$\|\alpha_n- \alpha\|_l\leq 
1/p_n\sum_{k=1}^{p_n}\|\varphi_n^k\,{}^*\alpha_n-\alpha\|_l\leq m_n,$$
where $m_n=\sup_{1\leq k\leq 
p_n}\|\varphi_n^k\,{}^*\alpha_n-\alpha\|_l$. For all $\nN$, there exists an integer $k_n$ such that $m_n=\|\varphi_n^{k_n}\,{}^*\alpha_n-\alpha\|_l$, therefore
$m_n\rightarrow_{n\rightarrow\infty} 0$ and for
all 
$l\in \mathbb Z$, 
$\|\alpha_n-\alpha\|_l\rightarrow_{n\rightarrow\infty} 0$
hence $\alpha_n\rightarrow_{n\rightarrow\infty} \alpha$, in a
$C^\infty$ way. According to the work of D. Ebin, for all $n$
sufficiently large, the isometry group of $\alpha_n$ is conjugated to a
subgroup of the isometry group of $\alpha$. That is a
contradiction. Consequently $h_\infty$ has at least one non trivial
isometry.
$\Box$\hop
{\bf Remarks.}
\begin{enumerate}
\item 
This result is still true for  manifolds  on which the limits of approximately isometric
systems (the metrics $h_\infty$ and $g_\infty$) always possess isometries. Of course the only known example is $\mathbf T^2$.
\item
It is possible to show (cf. \cite{these} examples p.48-49) that, on the torus, the
result is not valid for the $C^0$ topology or that we can have a diminution of the number of isometries for a $C^\infty$-convergent sequence of metrics.
But we found  more interesting to show what happens to this set
on  the compact quotients of $PSL(2,\mathbb R)$.
\end{enumerate}
\subsection{On compact quotients of $PSL(2,\mathbb R)$}
The goal of this section is to construct a sequence of metrics on a compact manifold with
non trivial isometry groups which converges ($C^\infty$) to a metric with
trivial isometry group. We are going to choose the metrics in the set of
left-invariant metrics of $PSL(2,\mathbb R)$ and look at the situation
on a compact quotient.
Consequently, we have to determine the
isometry groups of left invariant metrics on $PSL(2,\mathbb R)$. The
isometry groups of left invariant Riemannian metrics on compact
simple groups are well known since the works of  T. Ochiai 
and T. Takahashi (cf. \cite{octa}) and the works of J. E. D'Atri and W. Ziller (cf. \cite{Atzil}
theorem 5 page 24). Adapting some of their results to  $PSL(2,\mathbb R)$, 
we compute those groups, and thus the sequence of metrics will appear quite naturally.
\subsubsection{Isometry groups of left invariant metrics of
$PSL(2,\mathbb R)$}
During this paragraph, we will denote by $G_0$ the connected component
of the identity  of a group $G$. We first prove a $PSL(2,\mathbb R)$ version of
the theorem  of D'Atri and Ziller.
\begin{pro}\label{aut}
Let $h$ be a left invariant, non biinvariant, pseudo-Riemannian metric
on $G=PSL(2,\mathbb R)$. Then the isometry group of $h$ is generated by
the left translations and a subgroup of the right translations.
\end{pro}
{\bf Proof:} Let Is$(h)$  be the isometry group of $h$ and $L(G)$ the subgroup of
Is$(h)$ composed of the left translations.\\
We first study  the connected component of the identity of  Is$(h)$
that we denote by  Is${}_0(h)$. We denote by Int$(G)$ the group of
inner automorphisms of $G$. We are going to use the following result:
\begin{theo}[Ochiai and Takahashi (cf. \cite{octa} theorem $3$)]
 \quad Let  $G$ be a semisimple Lie group endowed with a  left invariant pseudo-Riemannian metric 
$h$. Then Is${}_0(h)\subset L(G)\,\mrm{Int}(G)$ if and only if 
 $L(G)$ is a  normal subgroup of Is${}_0(h)$.
\end{theo}
We are going to show that  $L(G)$ is a  normal subgroup of
Is${}_0(h)$. Let us give an upper bound for the dimension of Is$(h)$. 
It is well known that the dimension of the isometry group
of an $n$ dimensional pseudo-Riemannian manifold is at most
$n(n+1)/2$. Moreover we know that in this case the sectional
curvature is constant. In our case it means that the metric is
biinvariant. So, here, the dimension is at most $5$. Moreover, if the isometry group
contains a $5$ dimensional subgroup it also implies that the curvature is
constant and the metric is biinvariant (it can be proved  directly
from a study of the action of the isometry group on the 2-Grassmannian but it 
can also be
deduced from the systematic study of the isometry groups of Lorentz
$3$-manifolds done by C. Bona et B. Coll
in \cite{boco}). 
Therefore we have dim  Is${}_0(h)=3$ or $4$. We denote by  $\mathfrak I$
the Lie algebra of Is${}_0(h)$ and  $\mathfrak g$ the subalgebra of
$\mathfrak I$ corresponding to $L(G)$. We want to show that $\mathfrak g$ has to be  
an ideal of $\mathfrak I$. If dim $\mathfrak I=3$ it is obvious. Let us study the case where 
dim $\mathfrak I=4$.  The Lie algebra cannot be semisimple. Thus its radical, 
$\mrm{rad}(\mathfrak I)$,  is not trivial and its dimension is at most $1$ because it cannot
intersect $\mathfrak g$ which is simple. Let us take two elements $x$ and $y$ of $\mathfrak g$
and $u$ a non zero element of  $\mrm{rad}(\mathfrak I)$. According to Jacobi identity we have 
$$\left[[x,y],u\right] + \left[[u,x],y\right] +
\left[[y,u],x\right]=0.$$
However, since
$\mrm{rad}(\mathfrak I)$ is a one dimensional ideal, this implies 
$$\left[[u,x],y\right] + \left[[y,u],x\right]=0,$$
and therefore $\left[[x,y],u\right]=0$.
As $[\mathfrak g,\mathfrak g]=\mathfrak g$, we find that  $\mathfrak g$ is an ideal of 
 $\mathfrak I$ therefore Is${}_0(h)$ can be written $L(G)\times K$, with $K$ a one
 dimensional subgroup of the group of right translations whose Lie algebra is  
$\mrm{rad}(\mathfrak I)$.
\hop
To finish the proof, we use the following theorem of  J. E. D'Atri and 
W. Ziller. It is not exactly the same statement as the original version 
(compare with theorem 5 page 24 of \cite{Atzil}) but the proof is exactly the same.
\begin{theo}[D'Atri and  Ziller(cf. \cite{Atzil})]
Let $G$ be a simple Lie group on which any left invariant pseudo-Riemannian metric $h$ satisfies 
Is${}_0(h)\subset L(G)$\,Int$(G)$. Let	$\varphi$ be an isometry between two non biinvariant 
left invariant metrics $h$ and $h'$. Then $\varphi$ is an automorphism of $G$.
\end{theo}
This theorem implies that Is$(h)\subset \mrm{Aut}(G)$. To conclude we
 note  that
$\mrm{Int}(SL(2,\mathbb R))= \mrm{Aut}(SL(2,\mathbb R))$ 
(cf. \cite{Ov} for example).
$\Box$
\hop
{\bf Remark.} Thanks to this result finding the isometry group of a left invariant 
metric is now a problem of linear algebra. More precisely, we are going to see that finding the
 subgroup $K$ is the same as finding the intersection between the
 isometry groups of two quadratic forms of $\mathbb R^3$. Actually, we just want the  set of isometries that fix the neutral element of $PSL(2,\mathbb R)$). 
For a metric $h$ we will denote it by Isot$(h,e)$. 
Proposition \ref{aut} entails that if $h$ is a left invariant metric then Isot$(h,e)$ is contained in 
Int$\left(PSL(2,\mathbb R)\right)$. 
For the Killing metric Isot$(Killing,e)$ is the whole of $O(2,1)$ 
(see \cite{salein}, proposition 2.3.2.5, for a nice geometric study of the isometry group 
of the Killing form), moreover Isotr${}_0(Killing,e)$ (the connected component of the identity)
 is the whole group of inner automorphisms. Therefore in order to know  Is$(h)$ it is sufficient to 
compute $SO_0(2,1)\cap\,$Is$(h_e)$, where Is$(h_e)$ denotes the group of linear isometries of the quadratic
 form obtained by restricting $h$ to $T_eG$.
\\\\ \indent
If we  choose a left invariant trivialisation of the tangent bundle, we can represent
the Killing form by a symmetric matrix $K$ and a left invariant metric $h$ by a matrix
$H$. Hence we can associate to $h$ the endomorphism, $N_h$, of $\mathbb R^3$ whose matrix
 is $H\,K^{-1}$. It is well known that both quadratic forms
 are simultaneously diagonalizable if and only if $N_h$ is diagonalizable.\\
From this discussion, and some classical computations, we can deduce the following proposition
(see \cite{these} for the proof and a more detailed version).
\begin{pro}\label{theprop}
Let $h$ be a left invariant Lorentzian metric on $PSL(2,\mathbb R)$.
\begin{enumerate}
\item
If  all the 
eigenvalues of $N_h$ are distinct then Is$(h)$ is isomorphic to 
$PSL(2,\mathbb R)\times \mathbb Z_2$. 
\item
If $N_h$ is not diagonalizable and has two distinct eigenvalues or has
only one eigenvalue  and a 1-dimensional eigenspace then Is(h) is isomorphic to 
$PSL(2,\mathbb R)$.
\end{enumerate}
\end{pro}
The metrics corresponding to the second case possess a left invariant lightlike geodesic
foliation of codimension $1$ whose tangent is the  lightlike
plane that  the metric and the Killing form have in common. It is not hard to see 
(thanks to the precious proposition 3.18 of \cite{chebi}), that
left-invariant metrics which possess such lightlike $2$-foliations are
those which have a lightlike plane in common with the Killing form.
\subsubsection{On the compact quotients.}
On some compact manifolds, we want  to construct a sequence of
metrics $\sui{\overline h}$ with isometry group
$\mathbb Z_2$ which tends to a metric  $\overline h_\infty$ with trivial
isometry group. Of course, as we have seen in the proof
of theorem \ref{ouvert}, $\overline h_\infty$ will possess a lightlike
geodesic foliation of codimension $1$.  
This is quite similar to the situation described in proposition \ref{theprop}. 
Thus, we try to find our example on $PSL(2,\mathbb R)$ among left invariant metrics.\\
We choose a left invariant trivialisation of the tangent bundle of
$PSL(2,\mathbb R)$ such that the  matrix of the Killing form in it is:
$$
\begin{pmatrix}
0&1&0\\
1&0&0\\
0&0&1
\end{pmatrix}.
$$
For all $n>1$, we define $h_n$ as the left invariant metric whose matrix in
this trivialisation is 
$$\begin{pmatrix}
0      & \alpha & 1/n \\
\alpha & 0      & \gamma   \\  
1/n    & \gamma & \delta 
\end{pmatrix},$$
with $\alpha$, $\gamma$ et $\delta$ in $\mathbb R_+^*$.  It is easy
to see that, for all $n$, the matrix $N_{h_n}$ has three distinct
eigenvalues. Hence, according to proposition \ref{theprop}, their
isometry groups are all isomorphic to $PSL(2,\mathbb R)\times \mathbb
Z/2\mathbb Z$.
The sequence of quadratic forms defined by those matrices converges,
therefore the sequence of left invariant metrics converges
$C^\infty$. We denote the limit by $h_\infty$.
The matrix associated to $h_\infty$ is clearly:
$$
\begin{pmatrix}
0      &\alpha  & 0    \\
\alpha &0       & \gamma \\  
0      &\gamma  & \delta 
\end{pmatrix}.$$ 
We notice that  $N_{h_\infty}$ is not diagonalizable and has two
different eigenvalues, therefore,
according to proposition \ref{theprop}, the isometry group of $h_\infty$ is reduced
to the left translations. Hence, we have already lost an isometry. Let us take the quotient by a 
cocompact lattice of $PSL(2,\mathbb R)$ and see what happens.
\vskip 0.3cm
Let $\Gamma$ be a cocompact fuschian subgroup  of $PSL(2,\mathbb
R)$. We denote by $V$ the quotient manifold
$\Gamma\backslash$SL$(2,\mathbb R)$, by $\Sigma$ the compact hyperbolic
surface  $\Gamma\backslash\mathbb H^2$ endowed with its natural Riemannian
metric, by  $N(\Gamma)$ the normalizer of  $\Gamma$ in $PSL(2,\mathbb
R)$ and by $\overline{h}$ the metric induced on $V$ by a left
invariant metric $h$ of $PSL(2,\mathbb R)$.
Clearly if Is$(h)=SL(2,\mathbb R)\times K$, we have Is$(\overline
h)=N(\Gamma)/\Gamma\times K$. But we have a geometrical interpretation of
$N(\Gamma)/\Gamma$, it is the isometry group of $\Sigma$. However, in
genus greater than $2$, there exist (and it is the generic case)
hyperbolic surfaces with trivial isometry group (see for example \cite{buser}
in the proof of theorem 6.5.3).
If $\Gamma'$ is a group corresponding to such a surface it satisfies\,:
$$N(\Gamma')/\Gamma'=\{e\}.$$
Consequently, since in the case of  $h_\infty$ we have $K=\{\mrm{Id}\}$,
we have shown  that 
  $$\mrm{Is}(\overline h_\infty)=\{Id\}.$$
\\
The following statement is therefore proved. 
\begin{theo}\label{notnice}
 There exist compact manifolds $V$, of dimension greater than $2$ such that the set
of metrics without isometries is not an open subset of $\mca L(V)$ endowed with the $C^\infty$ 
topology.
\end{theo}
{\bf Remark.} Thanks to the article \cite{guelaf} of M. Guediri and J. Lafontaine, the same  sequence
 of  metrics enables us to see that the set of complete metrics is not closed in general
(it is, of course, not open).
It is not a surprise. Indeed, as we saw in lemma \ref{expo}, the result of 
\cite{rom-san} on the completeness of a geodesic obtained as a limit can 
be adapted to the completeness of a metric obtained as a limit.

\vskip .7cm\noindent
{\it
D\'epartement de Math\'ematiques, GTA (CNRS UMR 5030), Universit\'e
Montpellier~II, case 51, 34095 Montpellier, France.
\\
E-mail: {\tt mounoud@math.univ-montp2.fr}  }
\end{document}